\documentclass[12pt]{article}
\usepackage{amsmath}
\usepackage{amsfonts}
\usepackage{amsthm}
\usepackage{graphicx}
\graphicspath{{Figures/}} 
\usepackage{overpic}
\usepackage{amssymb} 
\usepackage{pictex}
\usepackage{rotating}  
\usepackage{color}
\usepackage{cite} 

\usepackage{enumitem} 

\usepackage{accents} 

\usepackage{yhmath} 
\usepackage{etex} 

\usepackage{comment}
\numberwithin{equation}{section}

\newtheorem{theorem}{Theorem}[section]
\newtheorem{proposition}[theorem]{Proposition}

\newtheorem{corollary}[theorem]{Corollary}
\newtheorem{Definition}[theorem]{Definition}

\newtheorem{Remark}[theorem]{Remark}

\newtheorem{RHproblem}[theorem]{RH problem}

\newtheorem{Example}[theorem]{Example}

\newcommand{\CC}{\mathbb{C}}

\newcommand{\RR}{\mathbb{R}}

\newcommand{\Z}{\mathbb{Z}}

\usepackage{color}

\renewcommand{\bar}{\overline}
\renewcommand{\tilde}{\widetilde}
\renewcommand{\hat}{\widehat}

\usepackage[bookmarksopen, naturalnames]{hyperref}
\begin{document}
\title{Correction/Addendum to ``The Extremal Function for the Complex Ball for Generalized Notions of Degree and Multivariate Polynomial Approximation''}
\author{T. Bloom, L. Bos, N. Levenberg, S. Ma'u and F. Piazzon}

\maketitle 

\section{Introduction} 
                  
                  In our Siciak memorial paper ``The Extremal Function for the Complex Ball ...'' \cite{paper} we construct the $P_{\infty}-$extremal function $V_{P_{\infty},B_2}$ for the complex Euclidean ball 
                  $$B_2:=\{(z_1,z_2): |z_1|^2+|z_2|^2\leq 1\} \subset \CC^2$$
                  and its complex Monge-Amp\`ere measure. Proposition 3.9 gives the formula for $V_{P_{\infty},B_2}$:
                  
 \begin{proposition}\label{ellinf} For $z=(z_1,z_2)$ with $|z_1|^2+|z_2|^2\ge1$,
\[V_{P_\infty,B_2}(z)=\begin{cases}
\frac{1}{2}\left\{\log(|z_2|^2)-\log(1-|z_1|^2)\right\}&\hbox{if  } |z_1|^2\le 1/2\,\,{\rm and}\,\,|z_2|^2\ge1/2 \cr
\frac{1}{2}\left\{\log(|z_1|^2)-\log(1-|z_2|^2)\right\}&\hbox{if  } |z_1|^2\ge 1/2\,\,{\rm and}\,\,|z_2|^2\le1/2 \cr
\log(|z_1|)+\log(|z_2|)+\log(2)&\hbox{if  } |z_1|^2\ge 1/2\,\,{\rm and}\,\,|z_2|^2\ge1/2 \end{cases}.
\]
\end{proposition}

We recall the notation and definitions in the next section. However our calculation of the Monge-Amp\`ere measure $(dd^c V_{P_{\infty},B_2})^2$, Proposition 4.2 in the paper, is incorrect. The error in our proof comes from the fact that the function $\frac{1}{2}[\log(|z_1|^2) - \log (1- |z_2|^2)]$ is {\it not} pluriharmonic in a neighborhood $U_a$ of any point  $a\in\partial B_2\cap \{|z_1|>1/\sqrt{2}>|z_2|\}$. In fact, it is fairly straightforward to see that the support of $(dd^c V_{P_{\infty},B_2})^2$ is the full topological boundary $\partial B_2$ of $B_2$ and the measure itself is absolutely continuous with respect to surface area measure on $\partial B_2$. Indeed, a much more general conclusion (Corollary \ref{correct}) can be obtained, giving a partial answer to the question of how the Monge-Amp\`ere measures $\mu_{P_q,B_2}:=(dd^c V_{P_q,B_2})^2$ vary with $q$ posed at the end of \cite{paper}. We give the correct calculation of $(dd^c V_{P_{\infty},B_2})^2$ in section 4. 

\section{Notation and definitions} Let 
\[P_1=\{(x_1,...,x_d)\in \RR^d: x_1,...,x_d \geq 0, \ x_1+\cdots x_d\leq 1\}\]
be the standard (unit) simplex in $(\RR^+)^d$. For a convex body $P\subset (\RR^+)^d$ we define the finite-dimensional polynomial spaces
\[Poly(nP):=\{p(z)=\sum_{J\in nP\cap (\Z^+)^d}c_J z^J: c_J \in \CC\}\]
for $n=1,2,...$ (here $J=(j_1,...,j_d)$); then $Poly(nP_1)$ are the usual polynomials of degree at most $n$. We consider $P$ satisfying (1.2) from \cite{paper}:
\begin{equation}\label{phyp} P_1 \subset kP \ \hbox{for some} \ k\in \Z^+.\end{equation}
We mention the particular examples 
$$P_q:=\{(x_1,...,x_d)\in (\RR^+)^d: (x_1^q+\cdots x_d^q)^{1/q}\leq 1\}, \ 1\leq q\leq \infty,$$ i.e., the (nonnegative) portion of an $l^q$ ball in $(\RR^+)^d$. 

The {\it indicator function} of a convex body $P$ is
\[\phi_P(x_1,...,x_d):=\sup_{(y_1,...,y_d)\in P}(x_1y_1+\cdots + x_dy_d).\]
For $P$ satisfying (\ref{phyp}) we have $\phi_P\geq 0$ on $(\RR^+)^d$ with $\phi_P(0)=0$.  
Define the logarithmic indicator function 
\[H_P(z):=\sup_{J\in P} \log |z^J|:=\phi_P(\log |z_1|,...,\log |z_d|)\]
where $|z^J|:=|z_1|^{j_1}\cdots |z_d|^{j_d}$ for $J=(j_1,...,j_d)\in P$. We use $H_P$ to define 
\[L_P=L_P(\CC^d):= \{u\in PSH(\CC^d): u(z)\leq H_P(z) +c_u \},\] (psh stands for plurisubharmonic) and 
\[L_P^+=L_P^+(\CC^d)=\{u\in L_P(\CC^d): u(z)\geq H_P(z) + C_u\}\]
where $c_u,C_u$ are constants depending on $u$. These are generalizations of the standard Lelong classes; i.e., when $P=P_1$. Given $E\subset \CC^d$, the {\it $P-$extremal function of $E$} is given by $V^*_{P,E}(z):=\limsup_{\zeta \to z}V_{P,E}(\zeta)$ where
\begin{equation}\label{vpk} V_{P,E}(z):=\sup \{u(z):u\in L_P(\CC^d), \ u\leq 0 \ \hbox{on} \ E\}.\end{equation}

For $K\subset \CC^d$ compact and nonpluripolar, we have
\begin{equation}\label{siczah} V_{P,K} =\lim_{n\to \infty} \frac{1}{n} \log \Phi_n \end{equation}
pointwise on $\CC^d$ where
\[\Phi_n(z):= \sup \{|p_n(z)|: p_n\in Poly(nP),  \ ||p_n||_K\leq 1\}\]
(cf., \cite{BHLevPer}).

From the definitions and from \cite{BBL}, we have the following for $K,K_1,K_2$ compact and $P,\tilde P$ satisfying (\ref{phyp}):
\begin{enumerate}
\item if $P\subset \tilde P$ then $V_{P,K}\leq V_{\tilde P,K}$;
\item if $K_1\subset K_2$ then $V_{P,K_1}\geq V_{P,K_2}$;
\item from (\ref{siczah}), $V_{P,K}=V_{P,\hat K}=V_{P,S_K}$ where 
$$\hat K:=\{z\in \CC^d: |p(z)|\leq ||p||_K, \ \hbox{all polynomials} \ p\}$$
is the polynomial hull of $K$ and $S_K$ is the Shilov boundary of $P(K)$, the uniform algebra generated by the polynomials restricted to $K$.
\end{enumerate}

Standard arguments from the definition (\ref{vpk}) imply that for $K$ compact and nonpluripolar, the support $S(P)$ of $\mu_{P,K}:=(dd^c V_{P,K}^*)^d$ is contained in $K$ (cf., \cite{BBL}). From 3. we have in fact $S(P)\subset S_K$. Finally, we say $K$ is {\it regular} if $V_{P_1,K}$ is continuous; this is equivalent to continuity of $V_{P,K}$ for all $P$ satisfying (\ref{phyp}) \cite{BBL}.

\section{Addendum} In the standard case $P=P_1$ (cf., \cite{BTfine}), for $K$ compact and nonpluripolar $S_K\setminus S(P_1)$ is pluripolar and is empty if $V_{P_1,K}$ is continuous. The same is true for general $P$ satisfying (\ref{phyp}).  
\medskip

\noindent {\bf Claim:} {\it The support $S=S(P)$ of $\mu_{P,K}:=(dd^c V_{P,K}^*)^d$ differs from $S_K$ by at most a pluripolar set and is independent of $P$}. 

\medskip
\noindent To verify this assertion, we begin with the following analogue of Proposition 5.3.3 \cite{K}. For $\Omega$ an open subset of $\CC^d$ and $E\subset \Omega$, the relative extremal function for $E$ relative to $\Omega$ is defined as \[\ u_{E,\Omega}(z)=\sup\{v(z): v\in PSH(\Omega),v|_E\leq -1, v\leq 0\}  \; , z\in \Omega. \]

\begin{proposition}
\label{pro:pro_3}

Let $P\subset (\RR^+)^d$ satisfying (\ref{phyp}), let $E$ be a bounded subset of $\mathbb{C}^d$, and let $\Omega$ be a bounded neighborhood of the polynomial convex hull of $\bar{E}$. If $E$ is not pluripolar, then there exist positive numbers $m$ and $M$ such that \[\ m(u_{E,\Omega}+1)\leq V_{P,E}\leq M(u_{E,\Omega}+1)\; \text{in}\; \Omega.\] Moreover, if $\Omega=\{V_{P,E}^*<C\}$ for some positive constant $C$, then \[\ V_{P,E}^*=C(u_{E,\Omega}^*+1)\; \text{in}\; \Omega. \]
\end{proposition}

\noindent The proof is straightforward and can be found in \cite{M}.

Theorem 7.1 of \cite{BTfine} shows, in particular, that for $K$ a compact, nonpluripolar, polynomially convex ($K=\hat K$) subset of a strictly pseudoconvex Runge domain $\Omega \subset \CC^d$ (e.g., a Euclidean ball), the support $S_{K,\Omega}$ of $(dd^cu^*_{K,\Omega})^d$ differs from $S_K$ by a pluripolar set and is empty if $u_{K,\Omega}$ is continuous. Thus in order to prove the assertion that $S_K\setminus S(P)$ is pluripolar and is empty if $V_{P,K}$ is continuous, it suffices to show $\mu_{P,K}:=(dd^c V_{P,K}^*)^d$ and $(dd^cu^*_{K,\Omega})^d$ are mutually absolutely continuous where $\Omega$ is any ball containing $\hat K$. This follows from Proposition \ref{pro:pro_3} and the following Monge-Amp\`ere comparison result which is a consequence of the proof of Lemma 2.1 in \cite{BT} (which itself is based on Lemma 2.1 in \cite{Lev}). The next proposition will also allow us to give some information on the variation of $\mu_{P_q,B_2}:=(dd^c V_{P_q,B_2})^2$ in $q$.

\begin{proposition} \label{lem21} Let $\Omega$ be a domain in $\CC^d$ and let $u,v\in PSH(\Omega)\cap L^{\infty}_{loc}(\Omega)$ satisfy the following properties:
\begin{enumerate}
\item there exists $S\subset \Omega$ compact and nonpluripolar with supp$(dd^cu)^d$=supp$(dd^cv)^d=S$ except perhaps a pluripolar set;
\item $S=\{u=0\}=\{v=0\}$ except perhaps a pluripolar set;
\item $u\geq v\geq 0$ in $\Omega$.
\end{enumerate}
Then $(dd^cu)^d\geq (dd^cv)^d$ as positive measures.
\end{proposition}

\noindent In \cite{BT} the set $S$ was contained in $\Omega \cap \RR^d$ but this hypothesis is not necessary (in that paper, they worked with subsets of $\RR^d$). We give a complete proof, following the arguments in \cite{BT}. Essentially, this is a generalization of Theorem 5.6.5 in \cite{K} where it is assumed $u,v$ are continuous and psh on a neighborhood of $\Omega$ and $u>v$ on $\partial \Omega$ (see also Lemma 2.1 of \cite{bkl}).

\begin{proof} We take standard smoothings $u_{\epsilon}=u*\chi_{\epsilon}$ and $v_{\epsilon}=v*\chi_{\epsilon}$. Let $\omega$ be a relatively compact subdomain of $\Omega$ with smooth boundary with $S\subset \omega$. Let $U,V$ be the restrictions of $u,v$ to $\partial \omega$ and let $U_{\epsilon},V_{\epsilon}$ be the restrictions of $u_{\epsilon}+2\epsilon,v_{\epsilon}+\epsilon$ to $\partial \omega$. Finally, for $\delta >0$, let 
$$S_{\delta}:=\{z: \hbox{dist}(z,S)\leq \delta\}.$$
We can assume $\delta$ is sufficiently small so that $S_{\delta} \subset \omega$. 

Define 
$$u_{\delta}^{\epsilon}(z):=\sup \{w(z): w\in PSH(\omega), \ w\leq 0 \ \hbox{on} \ S_{\delta}, \ \limsup_{z \to \zeta}w(z)\leq U_{\epsilon}(\zeta), \ \zeta \in \partial \omega \}$$
and 
$$v_{\delta}^{\epsilon}(z):=\sup \{w(z): w\in PSH(\omega), \ w\leq 0 \ \hbox{on} \ S_{\delta}, \ \limsup_{z \to \zeta}w(z)\leq V_{\epsilon}(\zeta), \ \zeta \in \partial \omega \}.$$
The set $S_{\delta}$ is a regular compact set and by the J. B. Walsh theorem (\cite{K}, Theorem 3.1.4), 
$u_{\delta}^{\epsilon}, v_{\delta}^{\epsilon}\in PSH(\omega)\cap C(\bar \omega)$. By Lemma 2.1 of \cite{Lev}, 
$$(dd^c u_{\delta}^{\epsilon})^d\geq (dd^c v_{\delta}^{\epsilon})^d.$$

Letting $\delta \downarrow 0$, we have 
$$u_{\delta}^{\epsilon} \uparrow u_{S,\epsilon}^* \ \hbox{a.e. in} \ \omega$$
where 
$$u_{S,\epsilon}(z):=\sup \{w(z): w\in PSH(\omega), \ w\leq 0 \ \hbox{on} \ S, \ \limsup_{z \to \zeta}w(z)\leq U_{\epsilon}(\zeta), \ \zeta \in \partial \omega \}.$$
Similarly,
$$v_{\delta}^{\epsilon} \uparrow v_{S,\epsilon}^* \ \hbox{a.e. in} \ \omega$$
where 
$$v_{S,\epsilon}(z):=\sup \{w(z): w\in PSH(\omega), \ w\leq 0 \ \hbox{on} \ S, \ \limsup_{z \to \zeta}w(z)\leq V_{\epsilon}(\zeta), \ \zeta \in \partial \omega \}.$$
From continuity of the Monge-Amp\`ere operator on locally bounded psh functions under a.e. increasing limits, we have
$$(dd^c u_{S,\epsilon}^*)^d\geq (dd^c v_{S,\epsilon}^*)^d.$$

Finally, we let $\epsilon \downarrow 0$. We have
$$u_{S,\epsilon}^* \downarrow u_S^* \ \hbox{on} \ \omega$$
where
$$u_S(z):=\sup \{w(z): w\in PSH(\omega), \ w\leq 0 \ \hbox{on} \ S, \ \limsup_{z \to \zeta}w(z)\leq U(\zeta), \ \zeta \in \partial \omega \};$$
and
$$v_{S,\epsilon}^* \downarrow v_S^* \ \hbox{on} \ \omega$$
where
$$v_S(z):=\sup \{w(z): w\in PSH(\omega), \ w\leq 0 \ \hbox{on} \ S, \ \limsup_{z \to \zeta}w(z)\leq V(\zeta), \ \zeta \in \partial \omega \}.$$
But clearly $u_S^*\leq u$ and $v_S^*\leq v$ on $\omega$; and since $u_S^*,v_S^*=0$ on $S$ except perhaps a pluripolar set, it follows that $u\leq u_S^*, \ v \leq v_S^*$ on $S$ except perhaps a pluripolar set. This set contains the supports of $(dd^cu_S^*)^d$ and $(dd^cv_S^*)^d$; hence by the domination principle (Corollary 4.5 \cite{btacta}) the reverse inequalities hold on $\omega$; i.e., $u_S^*\geq u$ and $v_S^*\geq v$ so that equality holds on $\omega$. By continuity of the Monge-Amp\`ere operator on locally bounded psh functions under decreasing limits, $(dd^cu)^d\geq (dd^cv)^d$.

\end{proof}

An immediate consequence of the Claim, together with Proposition \ref{lem21} and item 1. from the last section, is the following result.

\begin{corollary} Let $K\subset \CC^d$ be compact and nonpluripolar. Let $P\subset \tilde P$ be convex bodies in $(\RR^+)^d$ satisfying \ref{phyp}. Then 
$$\mu_{P,K}=(dd^c V_{P,K}^*)^d \leq \mu_{\tilde P,K}=(dd^c V_{\tilde P,K}^*)^d.$$
\end{corollary}

We now go back to the case of $K=B_2\subset \CC^2$. This is a regular compact set with $S_K=\partial B_2$. 

\begin{corollary}\label{correct} For $K=B_2\subset \CC^2$, the measures $\mu_{P_q,B_2}:=(dd^c V_{P_q,B_2})^2$ for $1\leq q \leq \infty$ all have support equal to the sphere $\partial B_2$. Moreover, for $q_1 \leq q_2$, 
$$\mu_{P_{q_1},B_2} \leq \mu_{P_{q_2},B_2}.$$
\end{corollary}

We know that $\mu_{P_{1},B_2}$ is a multiple of surface area measure $d\sigma$ on $\partial B_2$. To be precise, we are normalizing so that 
$$\int_{\CC^2} \mu_{P_{1},B_2}=\int_{\CC^2} (dd^c H_{P_1})^2=(2\pi)^2=4\pi^2$$
since we take $d=\partial +\bar \partial$ and $d^c =i(\bar \partial - \partial)$. Writing 
\begin{equation}\label{dsig} d\sigma = \bar z_1 dz_1\wedge d\bar z_2\wedge dz_2+\bar z_2 d\bar z_1 \wedge dz_1\wedge dz_2,\end{equation}
we have$$\int_{\partial B_2} d\sigma = 4 \pi^2 \ \hbox{so that} \ (dd^c H_{P_1})^2=d\sigma.$$
Here, in the first integral, using polar coordinates
$z_1=r_1e^{i\theta_1}$, $z_2=r_2e^{i\theta_2}$ with $\theta_1,\theta_2\in[0,2\pi]$, 
and  $r_1=\cos\psi$, $r_2=\sin\psi$ with $\psi\in[0,\pi/2]$ gives 
$$\int_{\partial B_2} d\sigma =\int_0^{2\pi}\int_0^{2\pi}\int_0^{\pi/2} 2\cos\psi\sin\psi d\psi d\theta_1 d\theta_2=4\pi^2.$$
Then for $1\leq q \leq \infty$, from \cite{BBL} we have
\begin{equation}\label{totalmass} \int_{\CC^2} \mu_{P_q,B_2} =\int_{\CC^2} (dd^c H_{P_q})^2 = 2! (4\pi^2) Vol(P_q)\end{equation}
where $Vol(P_q)$ denotes the Euclidean area of $P_q\subset (\RR^+)^2$.

\section{Calculation of $(dd^c V_{P_{\infty},B_2})^2$} To compute $(dd^c V_{P_{\infty},B_2})^2$, we use our explicit formula (\ref{ellinf}) and a result from \cite{BM}. From their main result, Theorem 1, defining
$$u(z_1,z_2)=\log(|z_2|^2)-\log(1-|z_1|^2),$$
it follows that there is a piece of $(dd^c V_{P_{\infty},B_2})^2$ supported on the portion of $\partial B_2$ where $|z_2|\geq |z_1|$ given by
$$\frac{1}{4}d^c u \wedge dd^c u.$$
We easily compute
$$\partial u = \frac{\bar z_1}{1-|z_1|^2}dz_1+\frac{1}{z_2}dz_2 \ \hbox{and} \ \bar \partial u = \frac{z_1}{1-|z_1|^2}d\bar z_1+\frac{1}{\bar z_2}d\bar z_2$$
so that
$$d^cu =i(\bar \partial u - \partial u)= i\bigl(\frac{z_1}{1-|z_1|^2}d\bar z_1+\frac{1}{\bar z_2}d\bar z_2-\frac{\bar z_1}{1-|z_1|^2}dz_1-\frac{1}{z_2}dz_2\bigr).$$
Then 
$$\frac{1}{i}dd^cu= d \frac{z_1}{1-|z_1|^2}\wedge d \bar z_1+ d\frac{1}{\bar z_2}\wedge d\bar z_2-d \frac{\bar z_1}{1-|z_1|^2}\wedge dz_1-d\frac{1}{z_2} \wedge dz_2$$
$$=\frac{1}{(1-|z_1|^2)^2}dz_1\wedge d \bar z_1-\frac{1}{(1-|z_1|^2)^2}d\bar z_1 \wedge dz_1=\frac{2}{(1-|z_1|^2)^2}dz_1 \wedge d\bar z_1.$$
Thus 
$$d^c u \wedge dd^c u =[\frac{z_1}{1-|z_1|^2}d\bar z_1+\frac{1}{\bar z_2}d\bar z_2-\frac{\bar z_1}{1-|z_1|^2}dz_1-\frac{1}{z_2}dz_2]\wedge \frac{2}{(1-|z_1|^2)^2}dz_1 \wedge d\bar z_1$$
$$=\frac{2}{z_2(1-|z_1|^2)^2}dz_2\wedge d\bar z_1 \wedge dz_1-\frac{2}{\bar z_2(1-|z_1|^2)^2}d\bar z_2\wedge d\bar z_1 \wedge dz_1.$$

To see that $d^c u \wedge dd^c u$ is a function multiple of $d\sigma$ in (\ref{dsig}), we first observe on $\partial B_2$ 
$$d(|z_1|^2+|z_2|^2) =0 = \bar z_1dz_1 +\bar z_2 dz_2+z_1d\bar z_1 + z_2 d \bar z_2$$ 
so that away from $z_1=0$,
$$d\bar z_1= \frac{-1}{z_1}[\bar z_1dz_1 +\bar z_2 dz_2+z_2d\bar z_2]$$
and
\begin{equation}\label{first} d\bar z_2\wedge d\bar z_1 \wedge dz_1= \frac{-\bar z_2}{z_1}d\bar z_2\wedge dz_2 \wedge dz_1.\end{equation}
Writing $d\sigma =d\sigma_1 +d\sigma_2$ where 
\begin{equation}\label{second}d\sigma_1 =\bar z_1 dz_1\wedge d\bar z_2\wedge dz_2 \ \hbox{and} \ d\sigma_2=\bar z_2 d\bar z_1 \wedge dz_1\wedge dz_2, \end{equation}
using (\ref{first}) we have
\begin{equation}\label{third}d\sigma_2=\frac{-|z_2|^2}{z_1}d\bar z_2 \wedge dz_1\wedge dz_2=\frac{|z_2|^2}{z_1}dz_1\wedge d\bar z_2\wedge dz_2=\frac{|z_2|^2}{|z_1|^2}d\sigma_1.\end{equation}
Hence for $|z_2|\ge |z_1|$,
\begin{equation}\label{fourth} d\sigma= (1+ \frac{|z_2|^2}{|z_1|^2})d\sigma_1.\end{equation}
Using (\ref{first}), then (\ref{second}) and (\ref{third}), and finally (\ref{fourth}), it follows that
$$d^c u \wedge dd^c u =\frac{2}{z_2(1-|z_1|^2)^2}dz_2\wedge d\bar z_1 \wedge dz_1 + \frac{2}{z_1(1-|z_1|^2)^2}d\bar z_2\wedge dz_2 \wedge dz_1$$
$$=\frac{2}{|z_1|^2(1-|z_1|^2)^2} \bar z_1 dz_1\wedge d\bar z_2\wedge dz_2 +\frac{2}{|z_2|^2(1-|z_1|^2)^2}\bar z_2 d\bar z_1 \wedge dz_1\wedge dz_2$$
$$=\frac{2}{(1-|z_1|^2)^2} \bigl[ \frac{1}{|z_1|^2}\bar z_1 dz_1\wedge d\bar z_2\wedge dz_2+ \frac{1}{|z_2|^2}\bar z_2 d\bar z_1\wedge dz_1\wedge dz_2\bigr].$$
$$=\frac{2}{(1-|z_1|^2)^2} \bigl[ \frac{1}{|z_1|^2}d\sigma_1+ \frac{1}{|z_1|^2}d\sigma_1\bigr]=\frac{4}{(1-|z_1|^2)^2}  \frac{1}{|z_1|^2}d\sigma_1$$
$$=\frac{\frac{4}{(1-|z_1|^2)^2}  \frac{1}{|z_1|^2}}{(1+ \frac{|z_2|^2}{|z_1|^2})}d\sigma=\frac{4}{(1-|z_1|^2)^2}  d\sigma=\frac{4}{|z_2|^4}d\sigma.$$

Similarly, it follows that there is a piece of $(dd^c V_{P_{\infty},B_2})^2$ supported on the portion of $\partial B_2$ where $|z_1|\geq |z_2|$ given by 
$$\frac{1}{4}d^c v \wedge dd^c v$$
where $v(z_1,z_2)=\log(|z_1|^2)-\log(1-|z_2|^2)$ so that
$$d^c v \wedge dd^c v=\frac{4}{|z_1|^4}d\sigma.$$
Thus the contribution to $(dd^c V_{P_{\infty},B_2})^2$ from these pieces is of the form $f_{\infty}(\zeta)d\sigma$ where $\zeta =(z_1,z_2)\in \partial B_2$ and
$$f_{\infty}(\zeta)=\frac{1}{\bigl(\max[|z_1|,|z_2|]\bigr)^4}.$$

Again using polar coordinates
$z_1=r_1e^{i\theta_1}$, $z_2=r_2e^{i\theta_2}$ with $\theta_1,\theta_2\in[0,2\pi]$, 
and  $r_1=\cos\psi$, $r_2=\sin\psi$ with $\psi\in[0,\pi/2]$, we compute
$$\int_{\partial B_2}f_{\infty}d\sigma= \int_{\partial B_2} \frac{1}{\bigl(\max[|z_1|,|z_2|]\bigr)^4}d\sigma$$
$$=\int_0^{2\pi}\int_0^{2\pi}\int_0^{\pi/4} \frac{ 2}{\cos^4 \psi}\cos\psi\sin\psi d\psi d\theta_1 d\theta_2+\int_0^{2\pi}\int_0^{2\pi}\int_{\pi/4}^{\pi/2} \frac{2}{\sin^4 \psi}\cos\psi\sin\psi d\psi d\theta_1 d\theta_2$$
$$=4 \int_0^{2\pi}\int_0^{2\pi}\int_0^{\pi/4} \frac{1}{\cos^4 \psi}\cos\psi\sin\psi d\psi d\theta_1 d\theta_2= 8\pi^2.$$
From (\ref{totalmass}), since $Vol(P_{\infty})=1$, this accounts for all of the mass of $(dd^c V_{P_{\infty},B_2})^2$. Thus 
$$(dd^c V_{P_{\infty},B_2})^2=\frac{1}{\bigl(\max[|z_1|,|z_2|]\bigr)^4}d\sigma.$$

\end{document}